\documentclass[12pt]{article}

\usepackage{geometry}                % See geometry.pdf to learn the layout options. There are lots.
\usepackage{amsmath,amsfonts,amssymb,amsthm,amscd,enumerate}
\usepackage[all]{xy}
\usepackage{graphicx}
\usepackage{epstopdf}

\usepackage[dvipsnames]{color}
\usepackage{relsize}

\title{Function algebras on the n-dimensional quantum complex space}

\author{{\sc Elmar Wagner} \\
\normalsize
Institute of Physics and Mathematics\\
\normalsize
University of Michoacan of San Nicolas of Hidalgo, Morelia, Mexico\\
\normalsize
e-mail: {\it elmar.wagner(at)umich.mx}\\
\mbox{ }\\
{\sc Ismael Cohen} \\
\normalsize
Universidad de la Costa, Barranquilla, Colombia\\
\normalsize
e-mail: {\it icohen(at)cuc.edu.co}}

\date{}                                           % Activate to display a given date or no date

\newtheorem{thm}{Theorem}%[section]
\newtheorem{prop}[thm]{Proposition}
\newtheorem{lem}[thm]{Lemma}

\theoremstyle{definition}
\newtheorem{defn}[thm]{Definition}

\newcommand{\nc}[2]{\newcommand{#1}{#2}}
\newcommand{\rnc}[2]{\renewcommand{#1}{#2}}
\rnc{\[}{\begin{equation}}
\rnc{\]}{\end{equation}}
\nc{\wegengruen}{\end{equation}}

\newcommand{\Z}{\mathbb{Z}}
\newcommand{\N}{\mathbb{N}}
\newcommand{\R}{\mathbb{R}}
\newcommand{\C}{\mathbb{C}}
\newcommand{\hsp}{{\hspace{-1pt}}}
\newcommand{\hs}{{\hspace{1pt}}}

\newcommand{\hH}{\mathcal{H}}

\newcommand{\G}{\mathcal{G}}
\newcommand{\cO}{\mathcal{O}}
\newcommand{\K}{{\mathcal{K}}}

\newcommand{\cL}{\mathcal{L}}

\newcommand{\CCq}{\cO(\C_q^2)} 
\newcommand{\Cqn}{\cO(\C_q^n)} 
\newcommand{\dd}{\mathrm{d}}

\newcommand{\im}{\mathrm{i}}
\newcommand{\e}{\mathrm{e}}

\newcommand{\spec}{\mathrm{spec}}
\newcommand{\dom}{\mathrm{dom}}

\newcommand{\supp}{\mathrm{supp}}

\newcommand{\rmB}{\mathrm{B}}

\newcommand{\lN}{\ell_2(\N)}

\newcommand{\ra}{\rightarrow}
\newcommand{\lra}{\longrightarrow}

\newcommand{\ip}[2]{\langle{#1},{#2}\rangle}
\newcommand{\msum}[2]{\underset{{#1}}{\overset{{#2}}{\mbox{$\sum$}}}}
\newcommand{\moplus}[2]{\underset{{#1}}{\overset{{#2}}{\mathlarger{\bar{\oplus}}}}}

\begin{document}
\maketitle

\begin{abstract}
The paper introduces a (universal) C*-algebra of continuous functions vanishing at infinity on the 
$n$-dimensional quantum complex plane. To this end, the well-behaved 
Hilbert space representations of the defining relations are classified. 
Then these representations are realized by multiplication operators on an $L_2$-space. 
The C*-algebra of continuous functions vanishing at infinity is defined 
by considering a *-algebra such that its classical counterpart separates the points of the 
$n$-dimensional complex plane and by taking the operator norm closure 
of a universal representation of this algebra. 
\end{abstract}
 
 \noindent 
{\it 2010 Mathematics Subject Classification:} Primary 46L65; secondary 58B32. 
\medskip
 
 \noindent 
{\it  Keywords and phrases:} n-dimensional quantum complex space, q-normal operators, 
C*-algebra generated by unbounded elements.

\section{Introduction}
Quantum groups, in their algebraic formulation, are generally defined as Hopf \mbox{(*-)al}\-gebras that are 
generated by a finite set of generators subjected to a number of algebraic relations \cite{KS}. 
In the same vain, the related quantum spaces, i.e.\ (*-)algebras on which a quantum group (co)acts, 
are then also given by a finite set of generators and relations. To study quantum groups and quantum spaces 
in a topological setting, it is desirable to associate a natural C*-algebra to the polynomial *-algebra 
of the generators. When the *-algebra admits faithful bounded representations, 
one may consider the closure of the polynomial *-algebra in the universal C*-norm given by the 
supremum of the operator norms of all (irreducible) Hilbert space representations. 
This strategy works well for compact quantum groups \cite{Wo} and compact quantum spaces such as 
Vaksman--Soibelman spheres and quantum projective spaces \cite{HS, VS}. 
%A prominent example is the Klimek--Lesniewski quantum disc, where the universal C*-al\-gebra 
%is isomorphic to the famous Toeplitz algebra \cite{KL}. 

The situation becomes more intricate if the polynomial *-algebra does not admit non-trivial 
bounded Hilbert space representations. To deal with such examples, Woronowicz introduced in \cite{Wo1} 
the concept of C*-algebras generated by unbounded elements. For non-compact  quantum groups 
or quantum spaces, these unbounded elements are usually provided by Hilbert space representations 
of the generators satisfying the required commutation relations on a dense domain. 
The C*-algebra generated by a set of closed operators (unbounded elements) 
has a natural interpretation as the $C_0$-functions on a locally compact quantum space. 
To get something like a universal C*-algebra, the unbounded elements should 
capture all possible Hilbert space representation, for instance by taking the direct sum  
or a direct integral of all irreducible representations. 

The obstacle of Woronowicz's approach is that it does not provide a constructive method, 
the C*-algebra has to be known beforehand in order to prove that it is generated by 
the unbounded elements. A proof of this property consists in showing that the unbounded 
generators are affiliated with the C*-algebra and provide information on all its 
Hilbert space representations, and that a product of the type 
$(1+t_1^*t_1)^{-1} \cdots (1+t_n^*t_n)^{-1}$ yields an element in the C*-algebra, where $t_j$ or $t_j^*$ 
belongs to the set of generators. This analysis has been performed in \cite{CW} for the 
1-dimensional complex quantum plane generated by a $q$-normal operator, 
i.e.\ a closed linear operator $z$ satisfying the operator equation $zz^*=q^2 z^*z$ 
for $q\in (0,1)$. Such operators have been studied systematically in \cite{CSS} and \cite{O}, in particular 
all their Hilbert space representations are known. Given a $q$-normal operator $z$, 
the generated C*-algebra can be described as a *-subalgebra 
of the crossed product algebra $C_0(\spec(|z|)) \rtimes\Z$, 
and this *-subalgebra can be defined in analogy to the classical case, see \cite{CW}.  
 
The 2-dimensional quantum complex plane has been studied by the authors in~\cite{CW1}. 
In this case, there are three types of Hilbert space representations: A trivial representation, which 
corresponds to the evaluation of functions at the origin; a representation where one of the generators is zero 
and the other is a $q$-normal operator, so this representation can be viewed as an embedding 
$\C_q\subset \C_q^2$; and finally a representation where the zero generator from the former representation 
is $q$-normal and injective, thus it may be viewed as describing functions 
on the Zariski open subset $\C_q^2\setminus \C_q$. The C*-algebra $C_0(\C_q^2)$ has then been defined 
by starting from the polar decomposition of the generators and describing a subalgebra of $C_0$-functions 
in analogy to the classical case such that classical counterpart separates the points of $\C^2$. 
Although all irreducible Hilbert space representation have been taken into account, 
a detailed analysis of Woronowicz's theory has not been carried out since we were only interested in 
a working definition for $C_0(\C_q^2)$. 

This paper treats the general case $C_0(\C_q^n)$ for any $n\in\N$. Again we are mostly interested 
in a working definition for $C_0(\C_q^n)$.  Nevertheless, in view of Woronowicz's approach, we will take all 
well-behaved Hilbert space representation into account. Roughly speaking, ``well-behaved'' means 
that we can apply the spectral theorem for certain distinguished self-adjoint operators and that 
the commutation relations hold in a strong sense. The well-behaved representations 
of the $n$-dimensional quantum complex plane are classified in Theorem \ref{reps}. 
Similarly to the case $n=2$, the different types of representations correspond to
embeddings $\{0\}\subset \C_q\subset \C_q^2 \subset \ldots \subset \C_q^n$. 
In Section \ref{fun}, we first realize these representations by multiplication operators 
on an $L_2$-space. To get something like a universal C*-algebra, we require that 
the measure has maximal support for each representation 
corresponding to an embedding $\C_q^k \subset \C_q^n$. 
Then we define a function *-algebra of $C_0$-functions on $\C_q^n$ such that 
its classical counterpart separates the points of $\C^n$. Finally the 
C*-algebra $C_0(\C_q^n)$ will be defined as the closure in the operator norm 
of a universal representation of this *-algebra.

\section{Preliminaries} 

Let $q\in(0,1)$ and $n\in\N$. By (the coordinate ring of) a $n$-dimensional quantum complex plane 
we mean the complex *-algebra $\Cqn$ 
generated by $z_1,\ldots,z_n$ satisfying the relations \cite{KS} 
\begin{eqnarray} \label{R1} 
&z_j z_i = q z_i z_j\hs, \quad j>i,   \ \qquad \ 
 z_j z_i^* = q z_i^* z_j \hs,   \quad j\neq i,  & \\   \label{R2}
&z_i z_i^* = q^2 z_i^* z_i - (1-q^2) \sum_{j>i}{}z_j^* z_j\hs, 
\quad i=1,...\hs,n\hsp -\hsp 1, 
\ \qquad\  z_n z_n^* = q^2 z_n^* z_n\hs.&
\end{eqnarray}

Our first aim will be to study *-representations of $\Cqn$, i.e., 
densely defined closed Hilbert space operators $z_1,\ldots, z_n$ and their adjoints $z_1^*,\ldots, z_n^*$ 
satisfying \eqref{R1} and \eqref{R2} on a 
common dense domain. A representation is called irreducible, if the Hilbert space does not decompose into the 
orthogonal sum of two non-trivial subspaces which are invariant under action of all elements from $\Cqn$. 
By a slight abuse of notation, we will use the same letter to denote 
a generator of the  coordinate ring $\Cqn$, its representation as a Hilbert space operator,  
and the restriction of this operator to a common dense domain. The so-called ``well-behaved'' representations will be 
classified in the next section, in this section we only provide preliminary results. 

To begin, note that the operator $z_n$ in \eqref{R2} satisfies the relation of a $q$-normal operator. Such operators have 
been studied in \cite{CSS,CW,OS,O,S}. In the next theorem, 
we restate the result from \cite[Corollary 2.2]{CW} with a slight change of 
notation: $n$ is replaced by $-n$ to obtain \eqref{zhn}.

\begin{prop} \label{corz}
Let $z$ be a non-zero $q$-normal operator, i.e., 
a densely defined closed linear operator 
on a separable Hilbert $\hH$ space satisfying the operator equation 
\[   												\label{zz*}
     z \hs z^* \,= \,q^2\hs z^*\hs z.
\]
Then the  Hilbert space $\hH$ decomposes into the direct sum 
$\hH=\ker({z})\oplus \underset{n\in\Z}{\bar\oplus} \hH_n$, where (up to unitary equivalence) $\hH_n = \hH_0$. 
The representation is uniquely determined by a self-adjoint operator $A$ on $\hH_0$ 
such that $\spec(A)\subset [q,1]$ and $q$ is not an eigenvalue. 
For $h\in \hH_0$, let $h_n$ denote the vector in $\hH$ which has 
$h$ in the $n$-th component of the direct sum $\underset{n\in\Z}{\bar\oplus} \hH_n$ 
and $0$ elsewhere. In this notation, the action of $z$ is determined by 
\[ \label{zhn}
z\hs h_n=  q^{-n} (A\hs h)_{n+1} 
\]
for all $h_n\in\hH_n$. The representation is irreducible, if and only if 
$\ker({z})=\{0\}$ and $\hH_0=\C$. In this case, $A$ stands for a real number belong to the interval $(q,1]$. 
\end{prop}

In the next lemma, the representations of another operator equation  are given 
that will be used in the classification of the well-behaved representation of $\Cqn$. 
The proof can be found in \cite{KW}. 

\begin{lem}                           \label{L1}
Let $w$ be a 
densely defined closed operator on a Hilbert space $\hH$
satisfying the operator equation 
\begin{equation}                                         \label{qhyp}
w w^\ast = q^2 w^\ast w  -(1-q^2). 
\end{equation} 
Then, up to unitary equivalence, $\hH$ decomposes into an orthogonal sum 
$\hH={\bar\oplus}_{n\in\N}\hH_n$, where $\hH_n =\hH_1$ for all $n\in \N$, 
and $w$ acts on $\hH$ by 
\[ \label{wh}
 w\hs h_n=  \sqrt{q^{-2n}-1}\hs h_{n+1}, \qquad h_n\in \hH_n. 
\]
The representation is irreducible if and only if $\hH_1=\C$. 
\end{lem}

The following simple lemma will be needed when the well-behaved representation are classified. 
\begin{lem}                           \label{L2}
 Given a *-representation of $\Cqn$ on a common dense domain $D$ of a Hilbert space $\hH$, let $\hH_0$ be a closed
 subspace of $\hH$ such that $D\cap \hH_0$ and $D\cap \hH_0^\bot$ are dense in $\hH_0$ and $\hH_0^\bot$, respectively. 
If $\hH_0$ is invariant under the action of $\Cqn$, it is reducing, i.e., $\hH_0^\bot$ is also invariant. 
\end{lem}
%
%\begin{proof}
The lemma follows from the fact that $a^* h \in \hH_0$ for all $h\in D\cap \hH_0$ and $a\in \Cqn$ 
since $\hH_0$ is invariant and $\Cqn$ 
is a *-algebra. Thus $\ip{ h}{a g} = \ip{ a^* h}{g} = 0$ for all $h\in D\cap \hH_0$ and $g\in D\cap \hH_0^\bot$, and hence  
$a g\in \hH_0^\bot$. 
%\end{proof}

%%%%%%%

%%%%%%%%

\section{Hilbert space representations}  \label{Hsr} 

In this section,  a complete description of the ``well-behaved" representations of $\Cqn$ will be given. 
It should be mentioned that no general method for defining well-behaved unbounded representations is known, 
so a case by case study is more common. The need of additional conditions to avoid pathological behavior 
has been shown convincingly in \cite[Sections 4.7 and 7.2]{S}. For $n=1$, the conditions on 
well-behaved representations of \eqref{zz*}  together with the related conclusions can be found in \cite{CSS} and \cite{CW}, 
as well as in \cite[Section 1.4.2]{OS} and \cite[Sections 11.3 and 11.6]{S}. 

To define what is meant by a ``well-behaved'' representations of $\Cqn$, 
we upgrade the analogous def\-i\-ni\-tions for $\CCq$ in \cite{CW1} to $\Cqn$. 
The main motivation of our definition is to be able to apply the spectral theorem for self-adjoint operators 
and to decompose the Hilbert space into the direct sum of (isomorphic) subspaces on which 
the operator relations from Proposition \ref{corz} and Lemma \ref{L1} can be realized.

As in \cite{CW1}, we start with some formal manipulations that motivate the definition. 
Hilbert space elements $h\in \hH$ are assumed to belong to the domains of the operators in consideration. 
For $k=1,\ldots,n$, set 
\[ \label{Qk}
Q_k := \msum{j=k}{n}\,z_j^* z_j\hs.
\]
Then obviously $Q_k^* = Q_k$, and 
\[ \label{zzQ} 
 z_k z_k^* - q^2z_k^* z_k = -(1-q^2)\hs Q_{k+1}
\]
by \eqref{R2}. Elementary calculations show that 
\begin{align} \label{zQ}
&&z_i Q_k&= Q_k z_i, & z_i^* Q_k&= Q_k z_i^*,  &  i&<k,&&\\ 
&&z_jQ_k&=q^2 Q_k z_j, & z_j^*Q_k&=q^{-2} Q_k z_j^*\,, & j&\geq k. \label{Qz}&&
\end{align}
Now, if $h\in\ker(Q_k)$, then $ Q_k (z_i h) =z_i (Q_k h)=0$ and $ Q_k (z_j h )= q^2 z_j (Q_k h)=0$, 
so $z_i h$ and $z_j h$ belong again to $\ker(Q_k)$. The same argument shows $z_i^* h,\, z_j^* h\in \ker(Q_k)$ for 
all $h\in\ker(Q_k)$. Assuming that the conditions of Lemma \ref{L2} are met, we can conclude that 
$\ker(Q_k)$ is reducing. 
On $\ker(Q_k)$, we have $z_k=\cdots = z_n=0$ by \eqref{Qk}. 
Setting $z_k=\cdots = z_n=0$ in \eqref{R1} and \eqref{R2}, one easily sees that the remaining generators 
fulfill the relations of $\cO(\C_q^{k-1})$. In particular, $z_{k-1}$ statisfies the relation \eqref{zz*} 
from Proposition \ref{corz}. Therefore it is natural to assume that $z_n$ on $\hH$ and $z_{k-1}$ on $\ker(Q_k)$, 
are $q$-normal operators. 

Later on,  we will classify the well-behaved representations by induction on $n$. Therefore we 
may assume that the well-behaved representations of $\cO(\C_q^{k-1})$ on $\ker(Q_k)$ are already known. 
Let us thus consider the representations on $\ker(Q_n)^\bot$, where $Q_1\geq\ldots \geq Q_n>0$. 
Here, $Q_k\geq Q_{k+1}$ follows from $Q_k = z_k^*z_k + Q_{k+1} $ by \eqref{Qk}. 
Equations \eqref{zQ} and \eqref{Qz} give for any polynomial $p$ in one variable 
\begin{align} \label{zpQ}
&&z_i p(Q_k) &= p(Q_k) z_i, & z_i^* p(Q_k) &= p(Q_k) z_i^*\hs, &  i&<k,&&\\ 
&&z_j p(Q_k) &=p(q^2 Q_k) z_j,  & z_j^*p(Q_k)&=p(q^{-2} Q_k) z_j^*\hs, & j&\geq k. \label{pQz}&&
\end{align}
Suppose that $Q_k$ is a self-adjoint operator and let $E_k$ denote the projection-valued measure such that 
$Q_k= \int \lambda\, \dd E_k(\lambda)$. Then $f(Q_k):= \int f(\lambda)\, \dd E_k(\lambda)$ yields for all 
measurable complex functions $f: \spec(Q_k) \ra \C$ a densely defined closed operator. 
We will assume that \eqref{zpQ} and \eqref{pQz} hold for  $p(Q_k)$ replaced by $f(Q_k)$
in an appropriate sense that takes care of the domains   
of unbounded operators, see Definition \ref{wb}.

On $\ker(Q_n)^\perp$, the relation $Q_k\hsp \geq\hsp Q_n\hsp >\hsp 0$ allows us to define 
 $\sqrt{Q_k}^{\hs -1} \hsp=\hsp \int\! \frac{1}{\sqrt{\lambda}} \,\dd E_k(\lambda)$. 
Consider for a moment 
\[
 w_k:= \sqrt{Q_{k+1}}^{\hs -1} z_k= z_k \sqrt{Q_{k+1}}^{\hs -1}    \label{w}
\]
as an abstractly defined operator. 
Inserting \eqref{w} into \eqref{zzQ} yields formally 
\[ \label{wwQ}
w_kw_k^* -q^2 w_k^*w_k = -(1-q^2) \sqrt{Q_{k+1}}^{\hs -1}  Q_{k+1} \sqrt{Q_{k+1}}^{\hs -1} = -(1-q^2). 
\]

Algebraically, \eqref{wwQ} is equivalent to \eqref{qhyp}. However, as in Lemma \ref{L2}, we
require that the operators $w_k$ satisfy \eqref{qhyp} as closed densly defined operators which includes 
the condition that $\dom(w_kw_k^*) = \dom(w_k^* w_k)$. 

Finally note that $z_1^* z_1,\ldots,z_n^* z_n\in\Cqn$ commute by \eqref{R1} and, as a consequence, 
so do $Q_1,\ldots,Q_n$ as elements in $\Cqn$. Viewed as self-adjoint operators, 
it is customary to require that they strongly commute. 

Similar to \cite[Definition 1]{CW1}, we sumarize our conditions on well-behaved representations 
in a definition.

\begin{defn} \label{wb}
A well-behaved *-representations of $\Cqn$ on a Hilbert space $\hH$ 
is given by densely defined closed operators $z_1,\ldots, z_n$ 
satisfying \eqref{R1}--\eqref{R2} on a common dense domain $D\subset \hH$ such that, for $j,k\in\{1,\ldots,n\}$, 
the following conditions hold: 
\begin{enumerate}[(i)]
\item $D$ is a core for the closed operators $z_j$ and $z_j^*$, 
and for the self-adjoint operators $Q_j$ given on $D$ by \eqref{Qk}. 

\item \label{0}
The self-adjoint operators $Q_j$ and $Q_k$ strongly commute. 

\item  \label{E} 
Let $E_j$ denote the unique projection-valued measure such that 
$Q_j= \int \lambda\, \dd E_j(\lambda)$. Then $E_j(M)\hs D \subset D$ 
for all Borel measurable sets $M\subset \R$. 

\item   \label{f} 
For all 
bounded Borel measurable functions $f$  on $\spec(Q_k)$, the operator relations 
\begin{align*}
&& f(Q_k) z_j&\subset z_j f(Q_k),  & f(Q_k) z_j^*&\subset z_j^* f(Q_k),   & j &< k,    &&   \\
&& f(Q_k) z_j&\subset z_j f(q^{-2}Q_j),  &  f(Q_k) z_j^*&\subset z_j^* f(q^2Q_k), & j&\geq  k, 
\end{align*} 
are satisfied. 

\item \label{zn} 
$z_n$ is a $q$-normal operator and, for $k<n$, the restriction of $z_k$ to $\ker(Q_{k+1})$ 
is a $q$-normal operator as defined in Proposition \ref{corz}. 

\item \label{Qzw}
On $\ker(Q_k)^\perp$, 
$z_j$ commutes with $\sqrt{Q_{j+1}}^{-1}$ for all $j=1,\ldots,k\hsp -\hsp 1$, and setting 
$w_j:= \sqrt{Q_{j+1}}^{-1} z_j= z_j \sqrt{Q_{j+1}}^{-1} $ determines a densely defined closed operator 
fulfilling the operator equation  
$$
w_jw_j^*  = q^2 w_j^*w_j -(1-q^2)  . 
$$ 
\end{enumerate}
\end{defn}

Definition \ref{wb}.\ref{E}) and \ref{f}) imply the following lemma.

\begin{lem} \label{red}
Let $k\in\{2, \ldots,n\}$ and $M_k,\ldots,M_n\subset [0,\infty)$ Borel measurable sets. 
Then for any  well-behaved *-representations of $\Cqn$ on  a common dense domain $D$ of 
a Hilbert space $\hH$, the subspace 
$E_k(M_k)\cdot\ldots\cdot E_n(M_n) \hH$ is reducing for the action of the operators 
$z_1,\ldots, z_{k-1}$ and $z_1^*,\ldots, z_{k-1}^*$, 
where  $E_j$ denotes the unique projection-valued measure from Definition \ref{wb}.\ref{E}). 
\end{lem}
\begin{proof}
Set $\hH_0:= E_k(M_k)\cdot\ldots\cdot E_n(M_n)\hH$. 
From Definition \ref{wb}.\ref{E}), we conclude that 
\begin{align*}
D\cap \hH_0\, &= \,E_k(M_k)\cdot...\cdot E_n(M_n)\hs D\, \subset \,D, \\
D\cap \hH_0^\bot\, &= \,\big(1-E_k(M_k)\cdot...\cdot E_n(M_n)\big)\hs D\,\subset \,D. 
\end{align*}
In particular, $D\cap \hH_0$ and $D\cap \hH_0^\bot$ are dense in $\hH_0$ and $\hH_0^\bot$, respectively. 
Moreover, Defi\-ni\-tion \ref{wb}.\ref{f}) implies $E_j(M_j) z_i \subset z_i E_j(M_j)$ since 
$E_j(M_j)  = \chi_{M_j}(Q_j)$, where $\chi_{M_j}$ stands for the indicator function of the set $M_j$. 
Hence, for all $h\in D\cap \hH_0$ and $g\in D\cap \hH_0^\bot$, and for all $i=1,\ldots,k-1$, 
we have 
\begin{align*}
& z_i\hs h = z_i\hs E_k(M_k)\cdot\ldots\cdot E_n(M_n)\hs h 
 = E_k(M_k)\cdot\ldots\cdot E_n(M_n)\hs z_i\hs h \in \hH_0, \\
 & z_i\hs g = z_i\hs \big(1-E_k(M_k)\cdot\ldots\cdot E_n(M_n)\big)g 
 = \big(1-E_k(M_k)\cdot\ldots\cdot E_n(M_n)\big) z_i\hs g \in \hH_0^\bot.
\end{align*}
Since the same holds for $z_i$ replaced by $z_i^*$, $\hH_0$ and $\hH_0^\bot$ are invariant 
under the densely defined actions of the operators 
$z_1,\ldots, z_{k-1}$ and $z_1^*,\ldots, z_{k-1}^*$, and therefore reducing. 
\end{proof} 

The next theorem is the main result of this section. It establishes a complete classification 
of the well-behaved *-representations of $\Cqn$. 

\begin{thm} \label{reps}
Any well-behaved *-representation of $\Cqn$ on a Hilbert space $\hH$ 
is unitarily equivalent to a representation 
described in the following way: For $n=1$, the well-behaved representations 
are outlined in Proposition 1. Now let $n>1$. Then the Hilbert space decomposes 
into the orthogonal sum of reducing subspaces 
\[ \label{hH}
\hH= \hH_0 \oplus \cdots \oplus \hH_n \ \ \ \text{such\ that}\ \ \ \ker(Q_k) = \hH_0 \oplus \cdots \oplus \hH_{k-1} , 
\]
where the self-adjoint operator $Q_k$ denotes the closure of \,$\sum_{j=k}^{n}\,z_j^* z_j$. 
Furthermore, there exist Hilbert spaces \,$\K_1,\ldots, \K_n$ 
such that each Hilbert space $\hH_k$, \,$k=1,\ldots,n$ decomposes into 
the orthogonal sum 
\[ \label{HK}
\hH_k = \moplus{i_1,...,i_{k-1}\in\N}{} \ \, \moplus{i_k\in\Z}{}\ \, \G_{i_1,...,i_{k-1},i_k}, \qquad 
\G_{i_1,...,i_{k-1},i_k}= \K_k. 
\]
The representation on $\hH_k$ is determined by a bounded self-adjoint operator $A_k$ on $\K_k$ 
such that $\spec(A_k) \subset [q,1]$, where $q$ is not an eigenvalue, and the actions of 
$z_1, \ldots, z_n$ on $\hH_k$ are given by 
\begin{align*}
 &z_{k+1} = \ldots =z_n= 0 \qquad  \text{if} \qquad k<n, \\
 &z_k h_{i_1,...,i_{k-1},i_k} = q^{-i_k}(A_k h)_{i_1,...,i_{k-1},i_k+1}\\
 &z_j h_{i_1,...,i_j, ... ,i_k} = \sqrt{q^{-2i_j}-1}\, q^{-(i_{j+1}+\cdots + i_k)}(A_k h)_{i_1,...,i_j+1, ... ,i_k},
 \qquad j=1,\ldots, k-1, 
\end{align*}
where $h\in \K_k$ and $h_{i_1,...,i_{k-1},i_k}:=h \in \G_{i_1,...,i_{k-1},i_k}=\K_k$. 
Moreover, 
$$
z_{1} = \ldots =z_n= 0 \qquad \text{on}\qquad \hH_0. 
$$ 
A common dense domain is given by $D:= D_0 \oplus \cdots \oplus D_n$, where $D_0:=\hH_0$ and, 
for $k>0$, %=1,\ldots,n$, 
\,$D_k:= \underset{i_1,...,i_{k-1}\in\N}{\oplus} \ \, \underset{i_k\in\Z}{\oplus}\ \, \G_{i_1,...,i_{k-1},i_k}$ 
is the algebraic orthogonal sum of $\G_{i_1,...,i_{k-1},i_k}= \K_k $. The representation is irreducible if all but 
one of the Hilbert spaces in the decomposition \eqref{hH} are zero and 
then either $\hH_0=\C$, or $\K_k=\C$ for the non-zero component $\hH_k$. In the latter case, $A_k\in (q,1]$ becomes a real number. 
\end{thm}
\begin{proof}
 We prove the theorem by induction on $n$. For $n=1$, Theorem \ref{reps} is equivalent to Proposition \ref{corz}, 
 and for $n=2$, the theorem was proven in \cite{CW1}. Now  let $n>1$ and assume that 
 the well-behaved representations of $\cO(\C_q^{n-1})$ are classified by Theorem~\ref{reps}. 
 As in Definition \ref{wb}.\ref{E}), we write 
 $Q_j= \int \lambda\, \dd E_j(\lambda)$ using the spectral theorem. 
 Then, by Lemma \ref{red}, $\ker(Q_n)= E_n(\{0\})\hH$ is reducing for 
 the operators $z_1,\ldots, z_{n-1}$ and $z_1^*,\ldots, z_{n-1}^*$. 
 Since  $\ker(Q_n)=\ker(z_n^*z_n)= \ker(z_n)$ 
 and, by \eqref{R2},   $\ker(z_n)= \ker(z_n^*z_n) =  \ker(z_n z_n^*)=\ker(z_n^*)$, we get $z_n=z_n^*=0$ on 
 $\ker(Q_n)$. In particular, $\ker(Q_n)$ is invariant under the action of all generators of $\Cqn$, 
 hence $\ker(Q_n)$ and $\ker(Q_n)^\bot$ are reducing by Lemma \ref{L2}. 
 Setting $z_n=z_n^*=0$ in the defining relations \eqref{R1} and \eqref{R2}, one sees that the remaining operators 
 satisfy the relations of $\cO(\C_q^{n-1})$. By induction, $\ker(Q_n)$ decomposes into the orthogonal sum 
 $\hH_0 \oplus \cdots \oplus \hH_{n-1}$, where the actions of $z_1,\ldots, z_{n-1}$  on $\ker(Q_n)$ are given in 
 Theorem~\ref{reps}. Finally, the Hilbert space $\hH$ decomposes into the orthogonal sum of reducing subspaces 
 $\hH= \hH_0 \oplus \cdots \oplus \hH_n$, where $\ker(Q_n)= \hH_0 \oplus \cdots \oplus \hH_{n-1}$ and 
 $\hH_n:= \ker(Q_n)^\bot$. Therefore it remains to classify the representations of $\Cqn$ on the reducing subspace 
 $\hH_n= \ker(Q_n)^\bot$.

 By Definition \ref{wb}.\ref{zn}), $z_n$ is a $q$-normal operator, therefore we can apply Proposition~\ref{corz}. 
 Since $\ker(z_n)= \ker(Q_n)$, we have $\ker(z_n)=0$ on $\hH_n$ and hence, by Proposition~\ref{corz}, $\hH_n$ 
 decomposes (up to unitary equivalence) into the orthogonal sum 
 $\hH_n= \underset{n\in\Z}{\bar\oplus} \G_{i_n}$, where $\G_{i_n} = \G_0$ for all $i_n\in \Z$. 
 Moreover, the actions of $z_n$ and $Q_n= z_n^* z_n$ are given by 
 $$
 z_n\hs h_{i_n}= q^{-i_n} (A_n h)_{i_n +1}, \qquad Q_n\hs h_{i_n}= q^{-2i_n} (A_n^2 h)_{i_n},\qquad h_{i_n}\in \G_{i_n}, 
 $$
 where $A_n$ denotes a bounded self-adjoint operator on $\G_0$ 
such that $\spec(A_n) \subset [q,1]$ and $q$ is not an eigenvalue. 
Note that $\spec(Q_n\!\!\upharpoonright_{\G_{i_n}})= q^{-2i_n}\spec(A_n^2)\subset [q^{-2i_n+2},q^{-2i_n}]$ with 
$q^{-2i_n+2}$ not being an eigenvalue. Writing $(0,\infty) = \cup_{i_n\in \Z}\, (q^{-2i_n+2},q^{-2i_n}]$ 
as the disjoint union of left-open and right-closed intervals, it can easily be seen that 
$$
\G_{i_n} = E_n((q^{-2i_n+2},q^{-2i_n}])\hH_n. 
$$

We proceed by induction, starting with the highest index $n$ and counting down.  
Let $j\in \{0,\ldots,n-2\}$ and suppose that 
$\hH_n = \moplus{i_{n-j},...,i_{n-1}\in\N}{} \ \, \moplus{i_n\in\Z}{}\ \, \G_{i_{n-j},...,i_{n-1},i_n}$
can be written as an ortho\-go\-nal sum of identical Hilbert spaces 
\,$\G_{i_{n-j},...,i_{n-1},i_n} = \G_{1,...,1,0}$. 
Assume furthermore that the actions of $z_{n-j},\ldots , z_n$ on above decomposition of $\hH_n$ are given by 
\[  \label{znhj}
 z_n h_{i_{n-j},...,i_{n-1},i_n} = q^{-i_n}(A_n h)_{i_{n-j},...,i_{n-1},i_n+1}
\]
and, if $j>0$, 
\[   \label{znlh}
z_{n-l} h_{i_{n-j},...,i_{n-l},..., i_n}  = \sqrt{q^{-2i_{n-l}}-1}\, q^{-(i_{n-l+1}+\cdots + i_n)}(A_n h)_{i_{n-j},...,i_{n-l}+1, ... ,i_n} 
\]
for $l=1,\ldots,j$. A direct computation shows that 
\begin{align*}
 & z_n^* z_n h_{i_{n-j},...,i_{n-1},i_n} = Q_n h_{i_{n-j},...,i_{n-1},i_n}   = q^{-2i_n}(A_n^2 h)_{i_{n-j},...,i_{n-1},i_n} ,\\
 &z_{n-l}^* z_{n-l} h_{i_{n-j},...,i_{n-l},...,i_n} 
 = (q^{-2i_{n-l}}-1)\, q^{-2(i_{n-l+1}+\cdots + i_n)}(A_n^2 h)_{i_{n-j},...,i_{n-l},...,i_n} . 
\end{align*}
Hence the action of $Q_{n-l} = \sum_{k=n-l}^n z_k^*z_k$ can be written as a telescoping sum yielding 
\[ \label{Qnl}
 Q_{n-l} h_{i_{n-j},...,i_{n-l},...,i_n}   = q^{-2(i_{n-l} + i_{n-l+1}+\cdots + i_n)}(A_n^2 h)_{i_{n-j},...,i_{n-l},...,i_n}. 
\]
Lastly, we obtain from \eqref{Qnl} that  
\begin{align} \label{EEH} 
\G_{i_{n-j},...,i_{n-1},i_n} =
E_{n-j} ((q^{-2(i_{n-j}+\cdots + i_n)+2},q^{-2(i_{n-j}+\cdots + i_n)}])\ \cdots\ E_n((q^{-2i_{n}+2},q^{-2i_n}])\hH_n .
\end{align}
This concludes the induction hypothesis. \medskip

To perform the induction step, we need to verify the following: \\[6pt]
1. $\G_{i_{n-j},...,i_{n-1},i_n}$ is reducing for the action of $z_{n-j-1}$. \\[6pt]
2. $\G_{i_{n-j},...,i_{n-1},i_n}$ is unitarily equivalent to the orthogonal direct sum 
\[ \label{Gnj}
 \G_{i_{n-j},...,i_{n-1},i_n} = \moplus{i_{n-j-1}\in\N}{}  \,\G_{i_{n-j-1},i_{n-j},...,i_{n-1},i_n}, 
 \quad \G_{i_{n-j-1},i_{n-j},...,i_{n-1},i_n} = \G_{1,1,...,1,0}, 
\]
3.  The actions of $z_{n-j-1}$ and $Q_{n-j-1}$ 
on $h_{i_{n-j-1},i_{n-j},...,i_{n-1},i_n}\in \G_{i_{n-j-1},i_{n-j},...,i_{n-1},i_n}$ read 
\begin{align} \nonumber
z_{n-j-1} h_{i_{n-j-1},i_{n-j},...,i_n}  
&= \sqrt{q^{-2i_{n-j-1}}-1}\, q^{-(i_{n-j}+\cdots + i_n)}(A_n h)_{i_{n-j-1}+1,i_{n-j},... ,i_n},\nonumber\\
 Q_{n-j-1} h_{i_{n-j-1},i_{n-j},...,i_n}  
&= q^{-2(i_{n-j-1}+i_{n-j}+\cdots + i_n)}(A_n^2 h)_{i_{n-j-1},i_{n-j},... ,i_n}. \label{3Q} 
\end{align}
4. Up to the unitary transformation mentioned in Item 2, we have 
\begin{align} \nonumber
 &\G_{i_{n-j-1},i_{n-j},...,i_n} 
 = E_{n-j-1} ((q^{-2(i_{n-j-1}+i_{n-j}+\cdots + i_n)+2},q^{-2(i_{n-j-1}+i_{n-j}+\cdots + i_n)}])\,\G_{i_{n-j},...,i_n} \\[2pt]
 &=E_{n-j-1} ((q^{-2(i_{n-j-1}+i_{n-j}+\cdots + i_n)+2},q^{-2(i_{n-j-1}+i_{n-j}+\cdots + i_n)}])  
 \nonumber  \\[2pt]
&\hspace{16pt} E_{n-j} ((q^{-2(i_{n-j}+\cdots + i_n)+2},q^{-2(i_{n-j}+\cdots + i_n)}])\ \cdots\ E_n((q^{-2i_{n}+2},q^{-2i_n}])\hH_n. 
\label{GEEH}
\end{align}
% For $l=0,\ldots,j$
5. The actions of $z_n,\ldots,z_{n-j}$ and $Q_n,\ldots,Q_{n-l}$ on 
$\G_{i_{n-j-1},i_{n-j},...,i_n} \subset \G_{i_{n-j},...,i_n}$ described in \eqref{znhj}--\eqref{Qnl} 
remain unchanged, 
the only difference being the appeareance of a new subindex $i_{n-j-1}$ corresponding to 
of Hilbert space decomposition \eqref{Gnj}.

Let  us first note that Item 1 follows immediately from \eqref{EEH} and Lemma \ref{red}. 
Since
$\sqrt{Q_{n-j}}^{-1}  E_{n-j} (M_j) = E_{n-j} (M_j) \sqrt{Q_{n-j}}^{-1}  E_{n-j} (M_j)$ 
by the spectral theorem, 
% for all Borel measurable sets $M\subset \R$, 
and since 
the bounded operator $\sqrt{Q_{n-j}}^{-1}  E_{n-j} (M_j)$ commutes with the 
projections $E_{n-l} (M_l)$ for all  $l=0,\ldots,j$ and all bounded, strictly positive 
Borel measurable sets $M_j,M_l\subset (0,\infty)$ by Definition \ref{wb}.\ref{0}), 
we conclude from Equation \eqref{EEH} and Item 1 that 
$w_{n-j-1}:= \sqrt{Q_{n-j}}^{-1} z_{n-j-1}$ also reduces $\G_{i_{n-j},...,i_n}$. 
Definition \ref{wb}.\ref{Qzw}) allows us to apply Lemma \ref{L1}. Therefore, up to unitary equivalence, 
$\G_{i_{n-j},...,i_n}$ decomposes into the orthogonal sum of identical Hilbert spaces 
\[   \label{G}
\G_{i_{n-j},...,i_n} = \moplus{k\in\N}{} (\G_{i_{n-j},...,i_n})_{k}, \qquad 
(\G_{i_{n-j},...,i_n})_{k} = (\G_{i_{n-j},...,i_n})_{1}\hs,
\]
and Equation \eqref{wh} gives
\[   \label{whj}
w_{n-j-1}(h_{i_{n-j},...,i_n})_{k} =  \sqrt{q^{-2k }-1}\,(h_{i_{n-j},...,i_n})_{k+1}, 
\]
where $(h_{i_{n-j},...,i_n})_{k} \in (\G_{i_{n-j},...,i_n})_{k}$. 

To deal with the operator $A_n$ on  $ (\G_{i_{n-j},...,i_n})_{k}$, consider the bounded measurable function  
$$
f: [0, \infty) \lra \R,\qquad 
f(t):=  q^{i_{n-j}+\cdots +i_n}\hs  \sqrt{t} \, 
\chi_{(q^{-2(i_{n-j}+\cdots +i_n)+2}\,,\,q^{-2(i_{n-j}+\cdots +i_n)}]}(t). 
$$ 
It follows from \eqref{Qnl} that $f(Q_{n-j}) = A_n$ on $\G_{i_{n-j},...,i_{n-1},i_n} = \G_{1,...,1,0}$. 
Furthermore, the restriction of $\sqrt{Q_{n-j}}^{-1}$ to $\G_{i_{n-j},...,i_{n-1},i_n}$ yields a bounded 
operator 
which obviously commutes with $f(Q_{n-j})$. From Definition \ref{wb}.\ref{f}), we get 
on $\G_{i_{n-j},...,i_{n-1},i_n}$ 
\begin{align} \nonumber
 A_n w_{n-j-1}^* w_{n-j-1}  &=f(Q_{n-j}) z_{n-j-1}^*Q_{n-j}^{-1} z_{n-j-1}  \subset 
     z_{n-j-1}^*Q_{n-j}^{-1} z_{n-j-1} f(Q_{n-j})\\
   &=   w_{n-j-1}^* w_{n-j-1} A_n . \label{Aw}
\end{align}
As well known, \eqref{Aw} implies that the eigenspaces $(\G_{i_{n-j},...,i_n})_{k}$ of 
$ w_{n-j-1}^* w_{n-j-1}$ are reducing for $A_n$ since, for all 
$(h_{i_{n-j},...,i_n})_{k} \in (\G_{i_{n-j},...,i_n})_{k}$, 
\begin{align*} 
w_{n-j-1}^* w_{n-j-1} A_n (h_{i_{n-j},...,i_n})_{k} 
&=A_n  w_{n-j-1}^* w_{n-j-1} (h_{i_{n-j},...,i_n})_{k} \\
&= (q^{-2k }-1) A_n (h_{i_{n-j},...,i_n})_{k} , 
\end{align*}
hence $A_n (h_{i_{n-j},...,i_n})_{k}$ belongs again to $(\G_{i_{n-j},...,i_n})_{k}$. 
Using   $A_n w_{n-j-1} \subset w_{n-j-1} A_n$, which can be proven as in \eqref{Aw}, 
we get from \eqref{whj} 
\begin{align}    \nonumber
A_n (h_{i_{n-j},...,i_n})_{k} &= \overset{k-1}{\underset{m=1}{\mbox{$\prod$}}}\mbox{$\frac{1}{\sqrt{q^{-2m }-1}}$} 
A_n w_{n-j-1}^{k-1}(h_{i_{n-j},...,i_n})_1 \\
 &= \overset{k-1}{\underset{m=1}{\mbox{$\prod$}}} \mbox{$\frac{1}{\sqrt{q^{-2m }-1}}$} 
 w_{n-j-1}^{k-1}A_n (h_{i_{n-j},...,i_n})_1 
 = (A_n (h_{i_{n-j},...,i_n})_1)_k    \label{An}
\end{align}
 for all $(h_{i_{n-j},...,i_n})_{k} \in (\G_{i_{n-j},...,i_n})_{k}$. 
Therefore the action of $A_n$ on $\G_{i_{n-j},...,i_n}$ is completely determined by its restriction 
to $(\G_{i_{n-j},...,i_n})_1$. 
By a slight abuse of notation, will denote this restriction again by $A_n$. 

Recall that all $\G_{i_{n-j},...,i_n}$ are identical, i.e., $\G_{i_{n-j},...,i_n}= G_{1,...,1,0}$. 
Changing in \eqref{G} the summation index $k$ to $i_{n-j-1}$ and defining 
$\G_{i_{n-j-1},i_{n-j},...,i_n} :=(\G_{i_{n-j},...,i_n})_{i_{n-j-1}}$, 
$h_{i_{n-j-1},i_{n-j},...,i_n}:=(h_{i_{n-j},...,i_n})_{i_{n-j-1}}$, we get 
$\G_{i_{n-j-1},i_{n-j},...,i_n} =(\G_{1,...,1,0})_{1}=\G_{1,1,...,1,0}$. 
In this notation, the action of  
$z_{n-j-1} = \sqrt{Q_{n-j}} w_{n-j-1}$ 
can be written 
\begin{align} \nonumber
 z_{n-j-1} h_{i_{n-j-1},i_{n-j},...,i_n} &= \sqrt{Q_{n-j}} w_{n-j-1}h_{i_{n-j-1},i_{n-j},...,i_n} \\
 &= \sqrt{q^{-2i_{n-j-1}}-1}\, q^{-(i_{n-j}+\cdots + i_n)}(A_n h)_{i_{n-j-1}+1,i_{n-j},... ,i_n},
\end{align}
where we used \eqref{Qnl} , \eqref{whj} and \eqref{An}. 
Furthermore, the same calculations which led to \eqref{Qnl} show that \eqref{3Q} holds. 
It is clear that taking the closure of $Q_{n-j-1}$ from \eqref{3Q} yields a self-adjoint operator. 
This finishes the steps 1)--3)  of the induction process. 

Now let $i_{n-j}, \ldots, i_n$ be fixed and $k\in\N$. From \eqref{3Q}, we conclude that 
the restriction of $Q_{n-j-1}$ to $\G_{k,i_{n-j},...,i_n}=(\G_{i_{n-j},...,i_n})_{k}$ 
is given by  $q^{-2(k+i_{n-j}+\cdots + i_n)}A_n^2$, in particular 
\begin{align} \nonumber
\spec(Q_{n-j-1}\!\!\upharpoonright_{\G_{k,i_{n-j},...,i_n}}) 
&\ = \ q^{-2(k+i_{n-j}+\cdots + i_n)}\,\spec(A_n^2) \\
&\ \subset \ [ q^{-2(k+i_{n-j}+\cdots + i_n)+2}\,,\,q^{-2(k+i_{n-j}+\cdots + i_n)}],  \nonumber
\end{align}
and $q^{-2(k+i_{n-j}+\cdots + i_n)+2}$ is not an eigenvalue. Therefore 
\[  \label{GinE}
\G_{k,i_{n-j},...,i_n} \subset 
E_{n-j-1}((q^{-2(k+i_{n-j}+\cdots + i_n)+2}\,,\,q^{-2(k+i_{n-j}+\cdots + i_n)}])\G_{i_{n-j},...,i_n},  
\]
where $E_{n-j-1}$ denotes the projection valued measure from Definition \ref{wb}.\ref{E}). 
Since $(q^{-2(k+i_{n-j}+\cdots + i_n)+2}\,,\,q^{-2(k+i_{n-j}+\cdots + i_n)}]
\cap (q^{-2(l+i_{n-j}+\cdots + i_n)+2}\,,\,q^{-2(l+i_{n-j}+\cdots + i_n)}] = \emptyset$
for all $k,l\in\N$ such that $k\neq l$, we have 
\begin{align} \label{GinG} 
 \G_{i_{n-j},...,i_n} &=\moplus{k\in\N}{} \G_{k,i_{n-j},...,i_n} \\
 &\subset 
\moplus{k\in\N}{} E_{n-j-1}((q^{-2(k+i_{n-j}+\cdots + i_n)+2}\,,\,q^{-2(k+i_{n-j}+\cdots + i_n)}])\G_{i_{n-j},...,i_n}  
 \subset 
\G_{i_{n-j},...,i_n}    \nonumber
\end{align}
From \eqref{GinE} and \eqref{GinG}, it follows that 
\[  \label{GisE}
\G_{k,i_{n-j},...,i_n} =
E_{n-j-1}((q^{-2(k+i_{n-j}+\cdots + i_n)+2}\,,\,q^{-2(k+i_{n-j}+\cdots + i_n)}])\G_{i_{n-j},...,i_n},  
\]
and inserting \eqref{EEH} into \eqref{GisE} gives \eqref{GEEH}. 

Finally, since $h_{i_{n-j-1},i_{n-j},...,i_n}  \in \G_{i_{n-j},...,i_n} $, 
the actions of the operators $z_{n-j}, \ldots,z_n$ and 
$Q_{n-j},\ldots, Q_n$ are given by \eqref{znhj}--\eqref{Qnl}, acting only on the indices $i_{n-j},\ldots,i_n$. 
This completes the induction step.

It remains to show that the representations from Theorem \ref{reps} are well-behaved. Direct calculations show 
that the operators satisfy the relations \eqref{R1} and \eqref{R2} on $D$. From the spectral theorem, it follws that 
$D$ is a core for $Q_j$, $|z_j|$ and $|z_j^*|$, and therefore also for $z_j$ and $z_j^*$, \,$j=1,\ldots,n$. 
Next, Definition \ref{wb}.\ref{E}) follows from \eqref{EEH}.  This, together with the fact the operator $A_n$ 
acting on $\G_{1,\ldots,1,0}$ is the same for all $Q_j$, implies Definition \ref{wb}.\ref{0}). 
From the formulas in Theorem \ref{reps}, it is obvious that the Hilbert spaces $\G_{i_1,...,i_n}$ are reducing 
for the actions of the self-adjoint operators $Q_k$, $|z_j|$ and $|z_j^*|$, and that the restrictions 
of these operators to $\G_{i_1,...,i_n}$ are bounded and commute. Using the spectral theorem, one easily sees that 
$\dom(|z_j^*|) = \dom(|z_j|)$ and 
$$
\dom(|z_j|) =  \{  \moplus{i_{1},...,i_{n-1}\in\N}{} \ \, \moplus{i_n\in\Z}{}\ \, h_{i_{1},...,i_n}  \in \hH_n     
: \msum{i_1,...,i_{n-1}\in \N}{} \,\, \msum{i_n\in\Z}{} \, \, \|\hs |z_j|\hs h_{i_1,...,i_n}\hs\|^2 \,<\, \infty \} .
$$
As $\|\hs |z_j|\hs f(Q_k) \hs h_{i_1,...,i_n}\hs\| \,=\, \|\hs f(Q_k) \hs |z_j|\hs h_{i_1,...,i_n}\hs\| 
\,\leq\, \|f\|_{\infty}\, \hs \|\hs |z_j|\hs h_{i_1,...,i_n}\hs\|$, we conclude that 
$\dom(f(Q_k) z_j) = \dom(|z_j| ) \subset \dom(|z_j|\hs  f(Q_k))=\dom(z_j f(Q_k))$ 
and analogously $\dom(f(Q_k) z_j^*) \subset \dom(z_j^* f(Q_k))$ for all $j, \,k\in\{1,\ldots,n\}$ and 
all bounded measurable functions $f$ on $[0,\infty)$. 
The relations from Definition \ref{wb}.\ref{f}) can now be checked directly. 
Similarly, using 
$$
 \|\hs \sqrt{Q_{j+1}}^{-1}   \hs z_j\hs h_{i_1,...,i_n}\hs\| \,=\,
\|\hs z_j\hs \sqrt{Q_{j+1}}^{-1}   \hs h_{i_1,...,i_n}\hs\| \,=\,  \|\hs |z_j|\hs \sqrt{Q_{j+1}}^{-1}   \hs h_{i_1,...,i_n}\hs\|, 
$$
we deduce that $\dom(\sqrt{Q_{j+1}}^{-1} z_j)= \dom(z_j\sqrt{Q_{j+1}}^{-1} )=\dom(|z_j|\hs \sqrt{Q_{j+1}}^{-1} )$
and that the operator equations in Definition \ref{wb}.\ref{Qzw}) hold. Finally, $z_n$ is a $q$-normal operator 
by Proposition \ref{corz}. 
\end{proof}

\section{C*-algebra of functions}  \label{fun} 

Our aim is to define a C*-algebra $C_0(\C_q^n)$ of continuous functions on the n-dimen\-sio\-nal
quantum complex plane vanishing at infinity. As in \cite{CW} and \cite{CW1}, the first step is to realize the 
representations on an $L_2$-space such that 
the self-adjoint  operators $Q_k$ and $|z_j|$ act as multiplication operators. 

\begin{thm} \label{rm}
The *-representations of $\Cqn$ in Theorem \ref{reps} are unitarily equivalent to a direct sum of representations of 
the following form: Set $I_q:= \{ q^{-m} : m\in\N\}$ and 
$X_{q,k}:= I_q^{k-1}\times [0,\infty)=I_q\times \ldots \times I_q\times [0,\infty)$. 
Let $\nu$ be a $q$-invariant Borel measure on $[0,\infty)$, i.e., 
$\nu(qM) = \nu(M)$ for all Borel measurable sets $M\subset[0,\infty)$.
Consider the product measure $\mu_k:= \delta \otimes \cdots \otimes\delta \otimes \nu$ 
on $X_{q,k}$, where $\delta$ denotes the counting measure on %the countable set 
$I_q$, and define $D_k:= \{ f\in \cL_2(X_{q,k},\mu_k) :  \supp(f)\subset I_q^{k-1}\times (0,\infty)\ \rm{is\ compact} \}$. 
For $k\in\{ 1,\ldots, n\}$, the operators $z_1,\ldots,z_n$ act on $h\in D_k$ by 
\begin{align}
 &z_{k+1} = \ldots =z_n= 0 \qquad  \text{if} \qquad k<n, \label{zkn0}\\
 &(z_k h)(t_1,\ldots,t_{k-1},t_k) = q t_k\hs h(t_1,\ldots,t_{k-1},qt_k),  \label{zkf} \\
 &(z_jh)(t_1,...,t_j, ... ,t_k) = \sqrt{(qt_j)^2-1}\, t_{j+1}\cdots  t_k\hs h(t_1,...,qt_j, ... ,t_k),
 \qquad j<k.   \label{zjf}
\end{align}
 \end{thm}
\begin{proof}
It suffices to assume that $k=n$, 
the other cases follow by setting $z_{k+1} = \ldots =z_n= 0$ and replacing $n$ by $k$. 
Given a representation of $\Cqn$ on $\hH_n$, we observe that $z_n$ acts on 
$$
\moplus{j\in\Z}{}\ \, \G_{i_1,...,i_{n-1},j} \cong \moplus{j\in\Z}{}\ \, \G_{1,...,1,j} =: \K 
$$
as a $q$-normal operator. It has been shown in  \cite{CSS} and \cite[Theorem 2.3]{CW} 
that the representation of $z_n = U\hs |z_n|$ on $\K$ is unitarily 
equivalent to a direct sum of representations on $\cL([0,\infty), \nu)\cong \K$ satisfying 
\[  \label{znf}
(|z_n| f)(t) = t f( t),\qquad (U f)(t) = f(q t),\qquad 
(z_n f)(t) =  qt f(q \hs t), 
\]
where $\nu$ denotes a $q$-invariant Borel measure on $[0,\infty)$. 
%, i.e., 
%$\nu(qM) = \nu(M)$ for all Borel measurable sets $M\subset[0,\infty)$. 
Note that 
the $q$-invariance of $\nu$ implies the  unitarity of $U$. 
Moreover, $\G_{1,...,1,j}$ may be taken to be $\cL_2((q^{j+1},q^{j}],\nu)$ 
and the isomorphism $\G_{1,...,1,j}\cong \G_{1,...,1,0}$ corresponds to the restriction 
of $U^{j}$ to that subspace, namely 
$U^{j}: \cL_2((q^{j+1},q^{j}],\nu) \ra \cL_2((q,1],\nu)$. 
%, where $U$ denotes the unitary operator 
%$$
%U : \cL([0,\infty), \nu) \lra  \cL([0,\infty), \nu), \quad (Uf)(t):= f(qt). 
%$$
%Note that the unitarity of $U$ follows from the $q$-invariance of $\nu$. 

Since $I_q$ is a countable discret set and $\delta$ is the counting measure, we get an explicit unitary isomorphism 
$\cL_2(I_q, \delta)\cong\lN$ by choosing the orthonormal basis $\{ \chi_m : m\in\N\}$ of $\cL_2(I_q, \delta)$, 
where $\chi_m:= \chi_{\{q^{-m}\}}$ denotes the indicator function of the singelton $\{q^{-m}\}$. 
The operator $S: \cL_2(I_q, \delta) \ra \cL_2(I_q, \delta)$ given by $(Sf)(t):= f(qt)$ becomes on this basis 
the shift operator, 
\[ \label{S} 
(S\chi_m)(t)= \chi_{\{q^{-m}\}}(qt) = \chi_{\{q^{-m-1}\}}(t) = \chi_{m+1}(t). 
\]
Furthermore, $\{ \chi_m : m\in\N\}$ is a basis of eigenvectors of the multiplication operator 
$T_\varphi: \cL_2(I_q, \delta) \ra \cL_2(I_q, \delta)$, $(Tf)(t) = \varphi(t)\hs f(t)$. More precisely, 
\[ \label{T} 
(T_\varphi\chi_m)(t)= \varphi(t)\hs \chi_{\{q^{-m}\}}(t) =\varphi(q^{-m}) \chi_{\{q^{-m}\}}(t) = \varphi(q^{-m})\chi_{m}(t), 
\]
since the support of $\chi_m$ is the singelton $\{ q^{-m}\}$. 
Next we have the following unitary equivalences of Hilbert spaces 
\begin{align}   \nonumber 
 \hH_n &=\! \moplus{i_1,...,i_{n-1}\in\N}{} \ \, \moplus{i_n\in\Z}{}\ \G_{i_1,...,i_{n-1},i_n} 
 =\! \moplus{i_1,...,i_{n-1}\in\N}{} \ \big( \moplus{j\in\Z}{}\ \, \G_{1,...,1,j} \big)
\cong \underbrace{ \lN \bar\otimes\cdots \bar\otimes  \lN}_{\substack{(n-1)\text{-times}}}\bar \otimes \K\\[-8pt]
&\cong  \cL_2(I_q, \delta) \bar\otimes\cdots \bar\otimes  \cL_2(I_q, \delta) \bar\otimes  \cL_2([0,\infty), \nu)  
\nonumber  \\
&\cong \cL_2(I_q\times \cdots \times I_q\times[0,\infty) , \delta\otimes \cdots \otimes \delta\otimes \nu) 
= : \cL_2(X_{q,n}, \mu_n), \label{ui}
\end{align}
where 
$$
\G_{i_1,...,i_{n-1},i_n}  \cong \chi_{i_1} \otimes \cdots \otimes \chi_{i_{n-1}} \otimes \cL_2((q^{i_n+1},q^{i_n}],\nu)
\cong  \cL_2((q,1],\nu), 
$$ 
and the last isomorphism in \eqref{ui} is obtained by multiplying the functions of the algebraic tensor product. 
Given a unitary isomorphism $\G_{1,...,1,0}\ni h\mapsto \hat{h}\in \cL_2((q,1],\nu)$ such that 
\[ \label{hat}
\widehat{(|z_n| h)}(t) = \widehat{A_n h}(t) =t\hat{h}(t), 
\]
as in \eqref{znf} and the paragraph following it, the unitary isomorphism \eqref{ui} is determined by 
\[ \label{hchi}
\G_{i_1,...,i_{n-1},i_n} \ni h_{i_1,...,i_{n-1},i_n} \ \longmapsto \ 
\chi_{i_1}  \cdots \chi_{i_{n-1}} U^{-i_n} \hat h =: \hat h_{i_1,...,i_{n-1},i_n} \in \cL_2(X_{q,n}, \mu_n)
\]
and the actions of $z_1,\ldots,z_n$ from \eqref{zkf} and \eqref{zjf} (with $k=n$) read on $\cL_2(X_{q,n}, \mu_n)$ as follows: 
\begin{align*}
 (z_n \hat h_{i_1,...,i_{n-1},i_n})(t_1,...\hs,t_{n-1},t_n) 
 &=q t_n  \chi_{i_1}(t_1)  \cdots \chi_{i_{n-1}}(t_{n-1}) ( U^{-i_n} \hat h)(q t_n) \\
 & = q^{i_n}  \chi_{i_1}(t_1)  \cdots \chi_{i_{n-1}}(t_{n-1}) \hs  U^{-i_n+1} (t_n \hat h( t_n)) \\
  & = q^{i_n}  \chi_{i_1}(t_1)  \cdots \chi_{i_{n-1}}(t_{n-1})\hs   U^{-i_n+1} ( \widehat{A_n h}( t_n)) \qquad\\
 &= q^{i_n}  ( \widehat{A_n h})_{i_1,...,i_{n-1},i_n-1}(t_1,...\hs,t_{n-1},t_n),  
 \end{align*}
 where we used \eqref{hchi} in the first and last equality, 
 \eqref{znf} in the second, and \eqref{hat} in the third. Analogously, applying additionally \eqref{S} and \eqref{T}, 
 \begin{align*}
 &(z_j \hat h_{i_1,...,i_{j},...,i_{n-1},i_n})(t_1,...\hs, t_j ,...\hs,t_{n-1},t_n) \\
 &=\sqrt{(qt_j)^2-1}\, t_{j+1}\cdots  t_n\hs  \chi_{i_1}(t_1) \cdots 
 \chi_{i_j}(qt_j) \cdots \chi_{i_{n-1}}(t_{n-1}) ( U^{-i_n} \hat h)(t_n) \\
  &=   \chi_{i_1}(t_1) \cdots 
 \big(\sqrt{(qt_j)^2-1}\,\chi_{i_j+1}(t_j)\big) \cdots  \big(t_{n-1}\chi_{i_{n-1}}(t_{n-1})\big) 
 \big(q^{i_n}\hs U^{-i_n} (t_n \hat h( t_n))\big) \\
 &=  \sqrt{q^{-2i_j}-1}\, q^{-(i_{j+1}+\cdots + i_{n-1}) +i_n} 
 \chi_{i_1}(t_1) \cdots \chi_{i_j+1}(t_j)\cdots \chi_{i_{n-1}}(t_{n-1})\hs   U^{-i_n} ( \widehat{A_n h}( t_n)) \\
 &= \sqrt{q^{-2i_j}-1}\, q^{-(i_{j+1}+\cdots + i_{n-1}) +i_n} 
 ( \widehat{A_n h})_{i_1,...,i_{j}+1,...,i_{n-1},i_n}(t_1,...\hs, t_j ,...\hs,t_{n-1},t_n). 
\end{align*}
Finally, replacing the last index $i_n$ by $-i_n$ shows that the representations given in \eqref{zkn0}--\eqref{zjf} 
are unitarily equivalent to the representations in Thereom \ref{reps}. 

To finish the proof, it remains to show that the domain $D_n$ from Theorem \ref{rm} is isomorphic to $D_n$ 
from Theorem \ref{reps}. This can easily by seen by observing that 
$$
\G_{i_1,...,i_{n-1},i_n}  \cong 
\{ f\!\!\upharpoonright_{\{q^{-i_1}\}\times\ldots \times\{q^{-i_{n-1}}\}\times (q^{i_n+1},\hs q^{i_n}]} \ : 
f \in \cL_2(X_{q,n},\mu_n) \}, 
$$
and $\supp\big(\, f\!\!\upharpoonright_{\{q^{-i_1}\}\times\ldots\times\{q^{-i_{n-1}}\}\times (q^{i_n+1},q^{i_n}]}\!\hsp\big)
\subset \{q^{-i_1}\}\times\ldots \times\{q^{-i_{n-1}}\}\times [q^{i_n+1},q^{i_n}]$ is compact. 
Moreover, if $f\in D_n\subset \cL_2(X_{q,n})$, its support is contained in a finite union of such sets. 
\end{proof} 

As in \cite{CW} and \cite{CW1},  the basic idea behind the definition of $C_0(\C_q^n)$ is to associate a 
*-algebra of bounded operators to well-behaved *-representations of $\Cqn$ 
such that its classical counterpart is dense in $C_0(\C^n)$. 
By taking direct sums and applying Theorem \ref{rm}, there is no loss of generality 
if we consider only representations on $\cL_2(X_{q,k},\mu_k)$ described in \eqref{zkn0}--\eqref{zjf}. 

Our starting point is the polar decomposition $z_j = S_j\hs |z_j|$, where $S_j$ and  $|z_j|$ act 
on $D_n\subset \cL_2(X_{q,n},\mu_n)$ by 
\begin{align}                                                                \label{SL2}
(S_jh)(t_1,...,t_j, ... ,t_n) &= h(t_1,...,q t_j, ... ,t_n),
 \qquad j=1,\ldots, n, \\
 (|z_n|\hs h)(t_1,\ldots,t_{n-1},t_n) &= t_n\hs h(t_1,\ldots,t_{n-1},t_n),  \label{znL2}\\ 
 (|z_j|\hs h)(t_1,...,t_j, ... ,t_n)                                           \label{zjL2}
 &= \chi_{{}_{I_q}\!}(t_j)\sqrt{t_j^2-1}\, t_{j+1}\cdots  t_n\hs h(t_1,...,t_j, ... ,t_n),\quad j<n, 
 \end{align}
where, although not necessary, we inserted the indicator function $\chi_{{}_{I_q}\!}$ 
such that the value $\chi_{{}_{I_q}\!}(q^k t_j)\sqrt{(q^kt_j)^2-1}$ remains well defined 
for all $k\in\Z$ and all $t_j \in I_q$. Obviously, $|z_1|, \ldots,|z_n|$ are commuting symmetric operators 
on $D_n$. For any polynomial $p\in \C[X_1,\ldots,X_n]$, the action of $p(|z_1|,...\hs,|z_n|)$ on $D_n$ 
is given by multiplication with the continuous function 
\[                                                                                \label{ptt}
\hat p(t_1,...\hs ,t_n) :=\, p\big(\chi_{{}_{I_q}\!}(t_1)\sqrt{t_1^2\!-\!1}\, t_{2}\cdots  t_n\,,\ldots, 
\chi_{{}_{I_q}\!}(t_{n-1})\sqrt{t_{n-1}^2\!-\!1}\, t_n\,, t_n\big)\hs \in\hs C( X_{q,n} ). 
\] 
As functions from $D_n$ have compact support, the action of $\hat p$ on $D_n$ is well defined 
and leaves $D_n$ invariant. Recall that $U$ in \eqref{znf} is a unitary operator, and $S$ in \eqref{S} 
acts as a unilateral shift with $\ker(S^*) = \chi_1= \chi_{\{q^{-1}\}}$. 
From these observations and Equation \eqref{SL2}, we get the following commutation relations on $D_n$:  
\begin{align}                                                          \label{Srel} 
& S_nS_n^* = S_n^* S_n =1, \qquad S_j^* S_j =1,\quad S_jS_j^* = 1-  \chi_{\{q^{-1}\}}(t_j), \quad j<n,  \\
& S_j S_k= S_k S_j, \qquad S_j S_k^*= S_k^* S_j, \qquad j\neq k, \\
&S_j\hs \hat p(t_1,...\hs ,t_j, ...\hs,t_n) = \hat p(t_1,...\hs ,qt_j, ...\hs,t_n)\hs S_j, \qquad 
 j=1,\ldots ,n.                                                     \label{Sp} 
\end{align}
The remaining relations are obtained by taking adjoints. 
Note that $1-  \chi_{\{q^{-1}\}}(t_j)$ belongs to $C( X_{q,n})$. 
For simplicity of notation, we set for $l\in\Z$
$$
S^{\# l}\,:=\,  \left\{   \begin{array}{cl} S^l\,,  & \ l \geq 0\,,   \\ 
                                      S^{*l}\,,  & \ l < 0\,.      \end{array}  \right.
$$
Then, using the operator relations \eqref{ptt}--\eqref{Sp}, one can easily show that the action of any polynomial 
$p(z_1,z_1^*,\ldots,z_n,z_n^*)$ on $D_n$ can be written as 
\[                                                                               \label{pSS}
p(z_1,z_1^*,\ldots,z_n,z_n^*) = \underset{\text{finite}}{\mbox{$\sum$}} \,
f_{l_1,...,l_n}(t_1,\ldots,t_n) S^{\# l_1}_1\cdots S^{\# l_n}_n,  \quad f_{l_1,...\hs,l_n}\in C( X_{q,n}). 
\] 

Equations \eqref{ptt} and \eqref{pSS} will be the main motivations for the definition of $C_0(\C_q^n)$. 
To explain this, let us take a look at the classical case. The polar decomposition $z_j = S_j \hs |z_j|$ 
corresponds to the Euler representation $z_j = \e^{\im \theta_j}\hs |z_j| $ of a complex number $z_j\in\C$. 
Given $l_1,\ldots,l_n\in\Z$ and $f_{l_1,...,l_n}\in C(\R_+^n)$, where $\R_+:=[0,\infty)$, 
the assign\-ment 
\[                                                                     \label{SW} 
\C^n \ni (\e^{\im \theta_1}\hs |z_1|,\ldots,\e^{\im \theta_n}\hs |z_n|)\ \longmapsto \ 
f_{l_1,...,l_n}(|z_1|,\ldots,|z_n|)\hs  \e^{\im l_1 \theta_1}\cdots \e^{\im l_n \theta_n} \in \C 
\]
defines an element in $C_0(\C^n)$ if and only if 
$f_{l_1,...,l_n}\in C_0(\R_+^n)$ and the function 
$f_{l_1,...,l_j,...,l_n}(|z_1|,...,0,...,|z_n|)\hs  
\e^{\im l_1 \theta_1}\hsp\cdot\!\cdot\!\cdot\hsp \e^{\im l_j \theta_j}\hsp\cdot\!\cdot\!\cdot\hsp \e^{\im l_n \theta_n}$ 
does not depend on $\theta_j$. % where the 0 occurs in the $j$-th position. 
Therefore we need to require 
\[           \label{fk0} 
f_{l_1,...,l_j,...,l_n}(|z_1|,...,0,...,|z_n|)=0 \quad \text{for} \quad l_j\neq 0. 
\]
On the other hand, by the Stone--Weierstra\ss\ theorem, the complex functions in \eqref{SW}, 
with $f_{l_1,...,l_n}\in C_0(  \R_+^n)$ satisfying 
the condition \eqref{fk0}, generate the C*-algebra $C_0(\C^n)$. For this reason, we will 
take analogous operators obtained from the representations in Theorem \ref{rm} 
as generators of $C_0(\C_q^n)$. 

Returning to the non-commutative case, we view the bounded operator $\pi_k(f)$ on 
$\cL_2(X_{q,k},\mu_k)$ given by multiplication with the function 
\begin{align} \label{pikf}
 \pi_k(f)&(t_1,...\hs, t_k,...\hs,t_n) \\ \nonumber  
 &:= f\big(\chi_{{}_{I_q}\!}(t_1)\sqrt{t_1^2\!-\!1}\, t_{2}\cdots  t_k\,,..., 
\chi_{{}_{I_q}\!}(t_{k-1})\sqrt{t_{k-1}^2\!-\!1}\, t_k\,, t_k,0,...\hs,0\big)
\end{align}
as an analog of $f(|z_1|,\ldots,|z_n|)\!\!\upharpoonright_{\R_+^k \times \{0\}}$ for $f\in C_0( \R_+^n)$, see \eqref{ptt}. 
For $k=0$, we define additionally $\pi_0(f)\in \rmB(\C)\cong \C\cong C(\{0\})$ by 
\[ \label{pikf0}
 \pi_0(f)(t_1,...\hs,t_n) := f(0,...\hs,0). 
\]
Furthermore,  for all $h\in \cL_2(X_{q,k},\mu_k)$, we set 
\[
 (\pi_k(S_j)h)(t_1,...,t_j, ... ,t_k) := h(t_1,...,q t_j, ... ,t_k), \ \ j\leq k, \quad \pi_k(S_j)=0,\ \  j>k, 
\]
and $\pi_k(S_j^*):= \pi_k(S_j)^*$.  
Finally, for all $l_1, \ldots,l_n\in\Z$ and 
$f_{l_1,...,l_n} \in C_0(  \R_+^n)$ satisfying \eqref{fk0}, we define 
\[   \label{pik}
 \pi_k(f_{l_1,...,l_n} S^{\# l_1}_1\cdots S^{\# l_n}_n) := \pi_k(f_{l_1,...,l_n})\pi_k(S_1)^{\# l_1}\cdots \pi_k(S_n)^{\# l_n}
\]
Note that, under the condition \eqref{fk0},  $\pi_k(S_j)=0$ for $j>k$ is consistent with \eqref{pikf} and \eqref{pikf0}. 
In this sense, the representation $\pi_k$ corresponds to the restriction of functions to $\C^k\times \{0\}\subset \C^n$. 

Next, recall that a C*-algebra generated by multiplication operators on a $L_2$-space 
does not depend on the measure, but only on the support of the measure. More precisey, 
assume that $\nu_j$ is a Borel measure on closed subsets $X_j\subset \R$  and let $M_j:=\supp(\nu_j)$. 
Then the C*-algebra generated by a family of multiplicacion operators 
$\{\pi(f_\alpha): \alpha\in J\}$, where 
$f_\alpha\in C_0(X_1\times \cdots\times X_k)$  and 
$$
(\pi(f_\alpha)h)(x_1,...\hs,x_k) := f_\alpha(x_1,...\hs,x_k) h(x_1,...\hs,x_k), \ \ 
h\in \cL(X_1\times ...\times X_k,\nu_1\otimes\hsp\cdot\!\cdot\!\cdot\hsp\otimes\nu_k), 
$$
is isomorphic to $C_0(M_1 \times ...\times M_k)$ if and only if the family 
$\{f_\alpha: \alpha\in J\}$ separates the points of $M_1 \times ...\times M_k$. 
Therefore, to get something like a universal C*-algebra (cf. \cite{CW} and \cite{CW1}), 
we will assume that the measure $\nu$ in Theorem \ref{rm} satisfies $\supp(\nu)=[0,\infty)$. 
 
Finally, recall that the representations of $\cO(\C_q^k)$ in Theorem \ref{reps} decompose into 
an orthogonal sum of representations on $\ker(|z_k|)$ and  $\ker(|z_k|)^\bot$. Moreover, 
the operator $|z_k|$ is invertible on $D_k\subset \ker(|z_k|)^\bot$. 
These representations can be viewed as representing functions on the disjoint sets 
$\{(z_1,...\hs ,z_k)\in\C^k_q : z_k=0\} = \C_q^{k-1}\times \{0\}$ and 
$\{(z_1,...\hs ,z_k)\in\C^k_q : z_k\neq 0\} =\C^k_q \setminus ( \C_q^{k-1}\times \{0\})$. 
Therefore, to obtain a C*-algebra that represents continuous functions on 
the whole quantum space $\C^n_q$, we will consider the orthogonal sum 
$ \pi_0\oplus \pi_1 \oplus ... \oplus \pi_n$. 
 
\begin{defn} 
Let $\nu$ be a $q$-invariant Borel measure on $[0,\infty)$ such that $\nu((q,1])\hsp<\hsp \infty$ and 
$\supp(\nu)= [0,\infty)$. For $k= 1,\ldots,n$, let $\hH_k:=\cL_2(X_{q,k},\mu_k)$ denote the 
Hilbert space given in Theorem \ref{rm} and set $\hH_0 =\C$. Then the 
C*-algebra $C_0(\C_q^n)$ of continuous functions on the n-dimen\-sio\-nal
quantum complex plane vanishing at infinity is the C*-subalgebra of 
$\rmB(\hH_0\oplus \ldots \oplus \hH_n)$ generated by the set of operators 
\begin{align*}
\big\{  (\pi_0\oplus \pi_1 \oplus ... \oplus \pi_n)\big(f_{l_1,...,l_n} S^{\# l_1}_1\cdots S^{\# l_n}_n\big) 
 :& \ l_1,\ldots,l_n \in \Z, \ \   f_{l_1,...,l_n}\in C_0(\R_+^n) ,  \\  
&\ f_{l_1,...,l_j,...,l_n}(r_1,\hsp...,0,\hsp...,r_n)=0\ \, \text{if}\, \ l_j\neq 0\hs\big\},
\end{align*} 
where $\pi_k(f_{l_1,...,l_n} S^{\# l_1}_1\cdots S^{\# l_n}_n)$ is defined in \eqref{pik}. 
\end{defn}

As in \cite{CW1}, we could have given an equivalent definition by considering a single representation on
$\cL_2(\R_+^n,\mu)$, where 
$$
\mu = \delta_0\otimes \cdots\otimes \delta_0 + \nu \otimes \delta_0 \otimes \cdots\otimes \delta_0 
+ \delta\otimes \nu \otimes \delta_0 \otimes \cdots\otimes \delta_0 + \ldots  
+\delta \otimes \cdots\otimes \delta\otimes \nu, 
$$
 $\delta$ denotes the extension of the counting measure on $I_q$ to the interval $[0,\infty)$, 
 and $\delta_0$ stands for the Dirac measure at the point 0. 
 However, the formulas for the representation would have been more cluttered involving 
 varios indicator functions in the arguments of the functions $f_{l_1,...,l_j,...,l_n}$, 
 see \cite[Definition 2]{CW1}. 
 
 The next aim will be to subject the C*-algebra $C_0(\C_q^n)$ to the same detailed analysis as 
 it has been done for $C_0(\C_q)$ in 
 \cite{CW} and \cite{CW2}, but that is behind the scope of this paper. 
 
 \section*{Acknowledgments} 
This work was partially supported by the 
CIC-UMSNH project ``Grupos cu\'anticos y geometr\'ia no conmutativa"
and by the CONACyT project A1-S-46784,  
and is part of the EU Staff Exchange project 101086394 “Operator Algebras That One Can See”.

\end{document}